\documentclass[10pt]{amsart}
\usepackage{palatino}
\usepackage[all]{xy,xypic}

\usepackage{amsfonts,amsmath,oldgerm,amssymb,amscd}

\newcommand{\ra}{\rightarrow}		
\newcommand{\lra}{\longrightarrow}
\newcommand{\by}[1]{\stackrel{#1}{\ra}}

\newcommand{\surj}{\ra\!\!\!\ra}	

\newcommand{\ol}{\overline}		
\newcommand{\wt}{\widetilde}
\newcommand{\iso}{\by \sim}

\newtheorem{theorem}{Theorem}[section]
\newtheorem{proposition}[theorem]{Proposition}
\newtheorem{lemma}[theorem]{Lemma}
\newtheorem{definition}[theorem]{Definition}
\newtheorem{corollary}[theorem]{Corollary}

\newcommand{\ga}{\alpha}

\newcommand{\bq}{\mbox{$\mathbb Q$}}

\newcommand{\ot}{\mbox{\,$\otimes$\,}}	
\newcommand{\op}{\mbox{$\oplus$}}
	
\newcommand{\hh}{\text{ht}}

\newcommand{\sur}{\twoheadrightarrow}
\newcommand{\bp}{\begin{proposition}}
\newcommand{\ep}{\end{proposition}}
\newcommand{\bl}{\begin{lemma}}
\newcommand{\el}{\end{lemma}}
\newcommand{\bt}{\begin{theorem}}
\newcommand{\et}{\end{theorem}}
\newcommand{\bc}{\begin{corollary}}
\newcommand{\ec}{\end{corollary}}
\newcommand{\bd}{\begin{definition}}
\newcommand{\ed}{\end{definition}}

\def\rmk{\refstepcounter{theorem}\paragraph{{\bf Remark} \thetheorem}}

\newcommand{\io}{\overline{I}}
\newcommand{\iw}{\widetilde{I}}

\newcommand{\wi}{\widetilde}

\def\rmk{\refstepcounter{theorem}\paragraph{{\bf Remark} \thetheorem}}

\def\proof{\paragraph{Proof}}
\def\example{\refstepcounter{theorem}\paragraph{{\bf Example} \thetheorem}}
\def\quest{\refstepcounter{theorem}\paragraph{{\bf Question} \thetheorem}}
\def\definition{\refstepcounter{theorem}\paragraph{{\bf Definition} \thetheorem}}
\newcommand{\remark}{\rmk}

\title{Efficient generation of ideals in a discrete Hodge algebra}
\author{Manoj K. Keshari}
\address{Department of Mathematics,
Indian Institute of Technology Bombay,
Powai, Mumbai 400076, India.}
\email{keshari@math.iitb.ac.in}
\author{Md. Ali Zinna}
\address{Department of Mathematics, Indian Institute of Technology Bombay, Powai, Mumbai 400076, India}
\email{zinna@math.iitb.ac.in}

\begin{document}
\maketitle

\begin{abstract}
Let $R$ be a commutative Noetherian ring and $D$ be a discrete Hodge
algebra over $R$ of dimension $d>\text{dim}(R)$. Then we show that 

(i) the top Euler class group $E^d(D)$ of $D$ is trivial.

(ii) if $d>\text{dim}(R)+1$, then $(d-1)$-st Euler class group
$E^{d-1}(D)$ of $D$ is trivial.
\end{abstract}

\section{Introduction}

Let $R$ be a commutative Noetherian ring. An $R$-algebra $D$ is called
a {\it discrete Hodge algebra over $R$} if $D=R[X_1,\cdots,X_n]/\mathcal{I}$,
where $\mathcal{I}$ is an ideal of $R[X_1,\cdots,X_n]$ generated by
monomials. Typical examples are $R[X_1,\cdots,X_n]$, $R[X, Y](XY)$
etc.  In \cite{v}, Vorst studied the behaviour of projective modules
over discrete Hodge algebras. He proved  \cite[Theorem 3.2]{v} that
{\it every finitely generated projective
  $D$-module is extended from $R$ if for all $k$, every finitely
  generated projective $R[X_1,\cdots,X_k]$-module is extended from
  $R$.  }

Later Mandal \cite{m2} and Wiemers \cite{wi} studied projective
modules over discrete Hodge algebra $D$. In \cite{wi}, Wiemers proved the
following significant result.
{\it 
Let $P$ be a projective $D$-module of 
rank $\geq \text{dim}(R)+1$. Then $(i)$ $P\simeq Q\oplus D$ for 
some $D$-module $Q$ and $(ii)$ $P$ is cancellative, i.e.
$P\oplus D\simeq P'\oplus D$ implies $P\simeq P'$.}

When $D=R[X,Y](XY)$, above results of Wiemers are due to Bhatwadekar
and Roy \cite{br}.  Very recent, inspired by results of
Bhatwadekar and Roy, Das and Zinna \cite{dz} studied the behaviour of
ideals in $R[X,Y]/(XY)$ and proved the following result on
efficient generation of ideals.
{\it Assume $\text{dim}(R)\geq 1$, $D=R[X,Y]/(XY) $
  and $I\subset D$ is an ideal of height $n=\text{dim}(D) $. 
Assume $I/I^2$ is generated by $n$ elements.
Then any given set of $n$ generators of $I/I^2$
  can be lifted to a set of $n$ generators of $I$.  In particular, the
  top Euler class group $E^n(D)$ of $D$ is trivial.}

As $R[X,Y]/(XY)$ is the simplest example of a discrete
Hodge algebra over $R$, motivated by above
discussions, one can ask the following question. \smallskip

\quest\label{ques0}
Let $R$ be a commutative Noetherian ring of dimension $\geq 1$ and 
$D$ be a discrete Hodge algebra over $R$ of dimension $n>\text{dim}(R)$.
 Let $I\subset D$ be an ideal of height $n$. Suppose that 
$I=(f_1,\cdots,f_n)+I^2$.
 Do there exist $g_1,\cdots,g_n\in I$ such that  $I=(g_1,\cdots,g_n)$ with $f_i-g_i\in I^2$?
 In other words, Is the top Euler class group $E^n(D)$ of $D$ trivial? (For definition of Euler class groups, see \cite{brs} and \cite{brs3}.) \smallskip

We answer Question \ref{ques0} affirmatively and prove the
following more general result ((\ref{main1}) below).  

\bp Let $R$ be a commutative Noetherian ring of dimension $\geq 1$ and $D$ be a discrete Hodge algebra over $R$ of dimension $n> \text{dim}(R)$.
 Let $P$ be a projective $D$-module of rank $n$ which is extended from $R$ and $I$ be an ideal in $D$ of height $\geq 2$. 
Suppose that there is a surjection $\alpha:P/IP\twoheadrightarrow I/I^2$. Then $\alpha$ can be lifted to a 
surjection $\beta: P \twoheadrightarrow I$.  In particular, the $n$-th Euler class group $E^n(D)$ of $D$ is trivial.
\ep

The above result can be extended to any rank $n$ projective $D$-module when $R$ contains $\mathbb{Q}$ (\ref{main0}). 
Here is the precise statement. 
\bt
Let $R$ be a commutative Noetherian ring containing $\mathbb{Q}$ of dimension  $\geq 2$ and $D$ be a discrete Hodge algebra over $R$ of dimension $n> \text{dim}(R)$.
 Let $I$ be an ideal in $D$ of height $\geq 3$ and $P$ be any rank $n$ projective $D$-module.
 Suppose that there is a surjection $\alpha:P/IP\twoheadrightarrow I/I^2$. Then $\alpha$ can be lifted to a 
surjection $\beta: P \twoheadrightarrow I$.
\et

After studying the top rank case, one is tempted to go one step further and 
inquire the following question. \smallskip

\quest\label{ques}
 Let $R$ be a commutative Noetherian ring of dimension $\geq 3$ and $D$ be a discrete Hodge algebra over $R$ of dimension
$d\geq \text{dim}(R)+2$. Let $I$ be an ideal in $D$ of height
$d-1$ and $P$ be a projective $D$-module of rank $d-1$.
 Suppose that $\alpha:P/IP \twoheadrightarrow
I/I^2$ is a surjection. Can $\alpha$ be lifted to a surjection
 $\beta:P\twoheadrightarrow I$? \smallskip

 We answer Question \ref{ques} affirmatively when $R$ contains $\mathbb{Q}$
 ((\ref{low}) below) as follows.  

 \bt\label{q1}
Let $R$ be a commutative Noetherian ring containing $\mathbb{Q}$ 
of dimension $\geq 3$ 
and $D$ be a discrete Hodge algebra over $R$ of dimension
$d>\text{dim}(R)$. Let $I$ be an ideal in $D$ of height
$\geq 4$ and $P$ be a projective $D$-module of rank $n\geq
max\{\text{dim}(R)+1, d-1\}$. Suppose that $\alpha:P/IP \twoheadrightarrow
I/I^2$ is a surjection. Then there exists a surjection
$\beta:P\twoheadrightarrow I$ which lifts $\alpha$. 
As a consequence, if $d\geq \text{dim}(R)+2$, then $(d-1)$-st Euler class 
group $E^{d-1}(D)$ of $D$ is trivial.
\et

Finally we derive an interesting consequence of above result as
follows (see (\ref{set})).

\bt
Let $R$ be a commutative Noetherian ring containing $\mathbb{Q}$ 
of dimension $\geq 3$  and $D$ be a discrete Hodge algebra over $R$ of
dimension $d>\text{dim}(R)$. Let $I$ be a locally 
complete intersection ideal in $D$ of height 
$n\geq max\{\text{dim}(R)+1, d-1\}$. 
Then $I$ is set theoretically generated by $n$ elements.
\et


\section{Preliminaries}

\noindent{\bf Assumptions.} Throughout this paper, rings are assumed to be commutative
  Noetherian and projective modules are finitely
  generated and of constant rank. For a ring A, 
$\text{dim}(A)$ will denote the Krull
  dimension of $A$. 
\smallskip

We start with the following definition.
\definition
An $R$-algebra $D$ is said to be a {\it discrete Hodge algebra over $R$} 
if $D$ is isomorphic to $R[X_1,\cdots,X_n]/J$, 
where $J$ is an ideal of $R[X_1,\cdots,X_n]$ generated by monomials.
A discrete Hodge algebra over $R$ is called {\it trivial} if it is a polynomial algebra over $R$. 
Otherwise, it is called a {\it non-trivial} discrete Hodge algebra.
\smallskip

\definition We call an ideal $I$ of a ring $R$ to be \emph{efficiently
  generated} if $\mu(I)=\mu(I/I^2)$, where $\mu(I)$
(resp. $\mu(I/I^2)$) stands for the minimal number of generators of
$I$ (resp. $I/I^2$) as an $R$-module (resp. $R/I$-module).

\smallskip

\definition Let $I$ be an ideal of a ring $R$. We say that $I$ is {\it
  set theoretically generated by $k$ elements} $f_1,\cdots,f_k$ in $R$
if $\sqrt {(f_1,\cdots,f_k)}=\sqrt I$.
\smallskip

The next two results are standard. For proofs the reader may consult \cite{brs}.

\bl\cite[2.11]{brs}\label{nak}
Let $R$ be a ring and $J$ be an ideal of $R$. Let $K\subset J$ and $L\subset J^2$ be two ideals of $R$ such that $K+L=J$. Then $J=K+(e)$ for some $e\in L$ with $e(1-e)\in K$ and $K=J\cap J^{\prime}$, where $J^{\prime}+L=R$. 
\el

\begin{lemma}\cite[2.13]{brs}\label{EE}
Let $A$ be a ring and $P$ be a projective $A$-module of rank $n$. Let 
$(\alpha, a)\in (P^\ast \oplus A)$. Then there exists an element 
$\beta \in P^\ast$ such that {\rm ht}\,$(I_a)\geq n$,
where $I = (\alpha + a\beta)(P)$. In particular, if the ideal 
$(\alpha(P),a)$ has height $\geq n$, then {\rm ht}\,$I\geq n$. Further, if 
$(\alpha(P),a)$ is an ideal of height $\geq n$ and $I$ is a proper ideal of 
$A$, then {\rm ht}\,$I=n$.  
\end{lemma}

The following lemma is proved in \cite[Lemma 3.1]{dk}.

\bl\label{move}
Let $R$ be a ring and $J\subset R$ be an ideal. Let $P$ be a projective $R$-module of rank 
$n\geq \text{dim}(R/J)+1$ and let ${\ga}:P/JP\sur J/J^2f$ be a surjection for some $f\in R$. 
Given any ideal $K\subset R$ with $\text{dim}(R/K) \leq n-1$, the map ${\ga}$ can be lifted to a surjection
 $\beta : P\sur J''$ such that: 
\begin{enumerate}
\item $J''+(J^2\cap K)f=J$,
\item $J''=J\cap J'$ and $\text{ht}(J')\geq n$,
\item $(J^2\cap K)f+J'=R$.
\end{enumerate}
\el

The following theorem is due to Mandal \cite[Theorem 2.1]{m}.
\bt \label{mandal}
Let $R$ be a ring and  $I\subset R[T]$ be an ideal containing a monic polynomial.
 Let $P$ be a projective $R$-module of rank $n\geq \text{dim}(R[T]/I)+2$. Suppose that there exists a surjection
$\phi:P[T]\twoheadrightarrow I/(I^2T).$
 Then, there exists a surjection $\psi:P[T]\twoheadrightarrow I$ which lifts $\phi$.
\et

We improve \cite[Lemma 2.9]{dz} in the following form to suit our needs. The proof is similar to the
one given in \cite[Lemma 4.9]{d}.

\bl\label{rem3.9} Let $R$ be a ring and $I$, $J$ be two ideals in $R$
such that $J\subset I^2$.  Let $P$ be a projective $R$-module and
$K\subset R$ be an ideal. Suppose that we are given surjections
$\alpha:P\twoheadrightarrow I/J$ and $\beta:P\twoheadrightarrow \ol I$
such that $\alpha\equiv\beta$ mod $\ol J$, where bar denotes reduction
modulo the ideal $K$.  Then $\alpha$ can be lifted to surjection
$\phi:P\twoheadrightarrow I/(JK)$.  \el
\medskip

The following result is implicit in the proof of  \cite[Theorem 3.2]{v}.

\bt\label{vorst}
Let $R$ be a ring and $r> 0$ be an integer. Assume that all projective modules of rank $r$ over polynomial extensions of $R$
 are extended from $R$. Then all projective modules of rank $r$ over discrete Hodge $R$-algebras  are extended from $R$.
\et
\smallskip

The following result is due to Das and Zinna \cite[Theorem 3.12]{dz1}.
\bt\label{subintegral}
Let $R$ be a ring of dimension $n\geq2$. 
 Let $R\hookrightarrow S$ be a subintegral extension and $L$ be a projective $R$-module of rank one. Then, the natural 
map  $E^n(R,L)\longrightarrow E^n(S,L\ot_R S)$ is an isomorphism.
\et

The following result follows from \cite[Lemma 3.2]{sw}.
\bl\label{sub}
Let $R\hookrightarrow S$ be a subintegral extension and $\mathcal{J}\subset R[X_1,\cdots,X_m]$ be an ideal generated by monomials. Then $R[X_1,\cdots,X_m]/\mathcal{J}\hookrightarrow S[X_1,\cdots,X_m]/\mathcal{J}$ is 
also subintegral. 
\el

The following result is from \cite[Proposition 2.13]{dz2}
for $d\geq 2$. By patching argument, it can be proved for
$d=1$.

\begin{proposition}\label{dz2.13}
Let $A$ be a ring of dimension $d \geq 1$. Let $I$ be an ideal of
$A[T]$ of height $\geq 2$ and $P$ be a projective $A[T]$-module of
rank $n \geq d +1$. Suppose that there exists a surjection $\phi : P/I P \surj
I/I^2$. Then $\phi$ can be lifted to a surjection $\Psi:P\surj I$.
\end{proposition}

The following result is due to Wiemers \cite[Corollary 4.3]{wi}.

\begin{theorem}\label{Wiemers}
Let $R$ be a ring of dimension $d$ and $D$ be a discrete Hodge algebra
over $R$. Let $P$ be a projective $D$-module of rank $>d$. Then

$(1)$ $P=D\oplus Q$ for some projective $D$-module $Q$.

$(2)$ $P$ is cancellative, i.e. if $P\oplus D\iso P'\oplus D$, then $P\iso P'$.
\end{theorem}

It is not hard to see that, adapting the same proof of \cite[Theorem
  4.2]{drs}, we can extend \cite[Theorem 4.2]{drs} in the following
form.

\bt\label{mrinal}
Let $R$ be a ring containing $\mathbb{Q}$ with $\text{dim}(R) = n\geq3$ and
$I\subseteq R[T]$ be an ideal of height $\geq 3$. Let $L$ be a projective $R$-module of rank $1$ and $P$ be a projective $R[T]$-module of rank $n$
whose determinant is $L[T]$. Assume that we are given a surjection
  $\psi: P\twoheadrightarrow I/(I^2T)$. Assume further that $\psi\otimes R(T)$ can be lifted to a surjection
 $\psi':P\otimes R(T)\twoheadrightarrow IR(T)$.
Then, there exists a surjection $\Psi: P \twoheadrightarrow I$ such that $\Psi$ is a lift of $\psi$ .
\et


\section{Main Theorems: Codimension Zero Case}

We begin with the following result which is motivated by
\cite[Theorem 4.2]{dz}.

\bp\label{main1}
Let $R$ be a ring of dimension $\geq 1$ and $D$ be a discrete Hodge algebra over $R$ of dimension $n> \text{dim}(R)$.
 Let $P$ be a projective $D$-module of rank $n$ which is extended from $R$ and $I$ be an ideal in $D$ of height $\geq 2$. 
Suppose that there is a surjection $\alpha:P/IP\twoheadrightarrow I/I^2$. Then $\alpha$ can be lifted to a 
surjection $\beta: P \twoheadrightarrow I$.
\ep

\proof If $D$ is a trivial discrete Hodge algebra over $R$, then we
are done by (\ref{dz2.13}). So we assume that $R$ is a non-trivial
discrete Hodge algebra over $R$. Let `prime' denote reduction modulo
the nil radical $N$ of $D$. Assume $\alpha\ot D'$ can be lifted to a
surjection $\alpha_1:P\ot D' \surj I\ot D'$. Then $\alpha_1$ can be
lifted to a surjection $\alpha_2:P_{1+N}\surj I_{1+N}$. Since $1+N$
consists of units of $D$, $\alpha_2$ is a lift of $\alpha$. Therefore,
we may assume that $D$ is reduced.

Let $D=R[X_1,\cdots,X_m]/{\mathcal{J}}$, where $\mathcal{J}$ is an
ideal of $R[X_1,\cdots,X_m]$ generated by square-free monomials. 
We
prove the result using induction on the number of variables $m$. If $m
= 1$, then $D$ is just $R[X_1]$ and the result follows from
(\ref{dz2.13}).

Let us assume that $m\geq 2$. We can assume that 
$\mathcal{J}=\mathcal{K}+X_m\mathcal{L}$, where $\mathcal{K}$ and
$\mathcal{L}$ are monomial ideals of $R[X_1,\cdots,X_{m-1}]$. Then
$D=R[X_1,\cdots,X_m]/(\mathcal{K},X_m\mathcal{L})$.

{\bf Case 1.} $n\geq 3$.
Given $\alpha:P/IP\twoheadrightarrow I/I^2$, applying (\ref{move}), 
$\alpha$ can be lifted to a
surjection $\gamma_1 : P\sur I'$ such that 
(1) $I'=I\cap J$,
(2) $I+J=D$,
(3) $\text{ht}(J)\geq n$.

If $\text{ht}(J)> n$, then $J=D$ and we are done.
So assume $\text{ht}(J)= n$. Let
$\gamma:P\surj J/J^2$ be the surjection induced from $\gamma_1$.

Let $x_m$ and $L$ be the images of $X_m$ and $\mathcal{L}$ in $D$, 
respectively. 
We shall use 'tilde' when we move modulo $(x_m)$ and 'bar' when 
we move modulo $L$.
  We first go modulo $x_m$ and 
consider the surjection $\wi \gamma:\wi P\twoheadrightarrow \wi J/\wi J^2$.
Note that $\widetilde J$ is an ideal of $\wt D=R[X_1,\cdots,X_{m-1}]/\mathcal{K}$ of height equal to dimension of 
 $\wt D$.
For this, we observe that $$\text{dim} (\wt D[X_m])+\text{ht}
(\widehat{X_m\mathcal L})=\text{dim}(D),$$ where 
$\widehat{X_m\mathcal L}$ is the image of $X_m\mathcal L$ in $\wt D[X_m]$.

 By induction hypothesis on $m$, there exists a surjection 
$\phi:\wi P\twoheadrightarrow \wi J$ which is a lift of $\wi\gamma$. 
Therefore, it follows from  (\ref{rem3.9}) that $\gamma$  can be lifted to a surjection $\psi:P\twoheadrightarrow J/(J^2x_m)$.

 We now move to the ring $\ol D=\dfrac
 {R[X_1,\cdots,X_{m-1}]}{(\mathcal{K},\mathcal{L})}[X_m]$ (i.e., go
 modulo $L$) and consider the surjection
$$\ol \psi:\ol P\twoheadrightarrow \ol J/(\ol J^2X_m)$$

Now observe that $J$ is of the form $J'/(X_m\mathcal{L})$ for some ideal 
$J'$ in $\frac{R[X_1,\cdots,X_{m-1}]}{\mathcal{K}}[X_m]$ containing 
$X_m\mathcal{L}$. 
Observe that $\text{ht}(J')= \text{dim}(\frac {R[X_1,\cdots,X_{m-1}]}{\mathcal{K}}[X_m])$. Therefore we may assume that $J'$ contains a monic polynomial 
in $X_m$. Since 
$\ol J=J/L\cap J= J'/L\cap J'$, it follows that $\ol J$ contains a monic in $X_m$. 
 Also $n\geq \text{dim}(\ol D/\ol J)+2(=2)$. 
By (\ref{mandal}), there exists a surjection $\theta:\ol P\twoheadrightarrow \ol J$ which lifts $\ol\psi$.

Therefore, it follows from (\ref{rem3.9}) that there exists a surjection
$\delta:P\twoheadrightarrow J/(J^2x_mL)$ which is a
lift of $\psi$.
As $x_m L=0$ in $D$, we obtain $\delta:P\twoheadrightarrow J$ is a surjection
which lifts $\gamma$. 
Now we have 
\begin{enumerate}
 \item $\gamma_1 : P\sur I\cap J$ such that $\gamma_1\ot D/I=\alpha\ot D/I$,
\item $\delta:P\twoheadrightarrow J$ with $\delta\ot D/J=\gamma_1 \ot
  D/J=\gamma$.
\end{enumerate}

Now by (\ref{Wiemers}), $P=D\oplus P'$. Also it follows that $n\geq \text{dim}(D/I)+2$
 and $n+\text{ht}(J)\geq \text{dim}(D)+3$. 
We can now use the subtraction principle \cite[Proposition 3.2]{dk} to find a surjection $\beta:P\twoheadrightarrow I$ 
 which lifts $\alpha$. This completes the proof in case $n\geq 3$.

{\bf Case 2.} $n=2$. In this case $\text{dim}(R)=1$ and hence by
(\ref{Wiemers}), $P\simeq L\op D$ for some rank one projective
$D$-module $L$.

We have $I=\alpha(P)+ I^2$. Applying (\ref{nak}), we can find $f \in I$ such that $I=(\alpha(P), f)$ with $f(1-f)\in \alpha(P)$ and 
 therefore we have a surjection $\alpha_{1-f}:P_{1-f}\twoheadrightarrow I_{1-f}$. 
 Let $\pi:P_f=L_f\op D_f\twoheadrightarrow D_f=I_f$ be the projection onto the second factor. Now consider the following surjections:
$$\alpha_{f(1-f)}:P_{f(1-f)}\twoheadrightarrow I_{f(1-f)}=D_{f(1-f)}$$
$$\pi_{1-f}:P_{f(1-f)}\twoheadrightarrow I_{f(1-f)}=D_{f(1-f)}$$ 

Now it is not hard to show that there exists $\tau\in SL(P_{f(1-f)})$ such that $\alpha_{f(1-f)}\tau=\pi_{1-f}$. Therefore standard patching argument
 implies that there is a projective $D$-module $Q$ of rank $2$ such that $Q$ maps onto $I$. By  (\ref{Wiemers}),
 $Q=\wedge^2(Q)\oplus D$. Also note that $Q$ has determinant $L$ and hence $Q\simeq L\op D$. 

By (\ref{Wiemers}), $L\op D$ is cancellative. We can now apply
\cite[Lemma 3.2]{b} to find a surjection $\beta:P\twoheadrightarrow I$
which lifts $\alpha$.  
\qed
\medskip



\bc\label{cor-main1}
Let $R$ be a ring of dimension $\geq 1$ and $D$ be a discrete Hodge algebra over $R$ of dimension $n> \text{dim}(R)$.
 Let $I$ be an ideal in $D$ of height $\geq 2$. Suppose that $I=(f_1,\cdots,f_n)+I^2$. Then there exist $g_1,\cdots,g_n$ such that 
$I=(g_1,\cdots,g_n)$ with $f_i-g_i\in I^2$ for $i=1,\cdots,n$.
\ec

\bc
Let $R$ be a ring of dimension $\geq1$ and $D$ be a discrete Hodge algebra over $R$ of dimension $n> \text{dim}(R)$. 
 Let $L$ be any rank one projective $D$-module.
 Then the $n$-th Euler class group $E^n(D, L)$ is trivial. 
\ec
\proof
 Let $D=R[X_1,\cdots,X_m]/\mathcal{J}$. 
Without loss of generality we can assume that $D$ is reduced (see \cite[Corollary 4.6]{brs}). In particular, $R$ is reduced. Let $S$ be the seminormalization of $R$ in its total 
quotient ring. Since $S$ is seminormal, by \cite[Theorem 6.1]{sw}, every rank one projective $S[X_1,\cdots,X_k]$-module is extended from $S$ for all $k$. 
Therefore, it follows from (\ref{vorst}) that $L\ot_R S$ is extended from $S$.

 Let us denote $S[X_1,\cdots,X_m]/\mathcal{J}$ by $D_1$. Since $R\hookrightarrow S$ is a subintegral extension, by (\ref{sub}),
 $D\hookrightarrow D_1$ 
is also subintegral. As $L\ot_R S$ is extended from $S$, by (\ref{main1}), it follows that $E^n(D_1, L\ot_R S)$ is trivial.
 Finally, using (\ref{subintegral}), we 
have $E^n(D,L)$ is trivial. 
\qed


\medskip

The following result is due to Katz \cite{ka}.
  
\bt\label{katz}
Let $R$ be a ring and $I\subset R$ be an ideal. Let $d$ be the
maximum of the heights of maximal ideals containing $I$, and suppose that $d < \infty$.
Then some power of $I$ admits a reduction $J$ satisfying $\mu(J/J^2)\leq d$.
\et

A result of Mandal from \cite{m2}, can now be deduced. 

\bc\label{sci}
Let $R$ be a ring of dimension $\geq 1$  and $D$ be a discrete Hodge algebra over $R$ of dimension $n> \text{dim}(R)$.
 Let $I\subset D$ be an ideal of height $\geq2$. Then $I$ is set theoretically generated by $n$ elements. 
\ec

\proof Using Katz (\ref{katz}), there exists $k>0$ such that $I^k$ has
a reduction $J$ with $\mu(J/J^2)\leq n$.  If $\mu(J/J^2)\leq
n-1$, then clearly $J$ is generated by at most $n$ elements. Therefore
we assume that $\mu(J/J^2)= n$. Since $J$ is a reduction of $I^k$, it
is easy to see that $\sqrt{I}=\sqrt{I^k}=\sqrt{J}$ and
$\text{ht}(I)=\text{ht}(J)$. Applying (\ref{cor-main1}), we see
that $J$ is generated by $n$ elements. Therefore, $I$ is
set-theoretically generated by $n$ elements.  \qed
\medskip

We have the following variant of (\ref{main1}) for rings containing $\bq$.

\bp\label{rem1}
Let $R$ be a ring containing $\mathbb{Q}$ of dimension  $\geq 2$ and $D$ be a discrete Hodge algebra over $R$ of dimension $n> \text{dim}(R)$.
 Let $I$ be an ideal in $D$ of height $\geq 3$ and $P$ be any rank $n$ projective $D$-module whose determinant is extended from $R$.
 Suppose that there is a surjection $\alpha:P/IP\twoheadrightarrow I/I^2$. Then $\alpha$ can be lifted to a 
surjection $\beta: P \twoheadrightarrow I$. 
\ep

\begin{proof}
We follow the proof of (\ref{main1}). The only thing which we need to
show is that $\ol \psi: \ol P\twoheadrightarrow \ol J/(\ol J^2 X_m) $
can be lifted to a surjection $\theta: \ol P\twoheadrightarrow \ol
J$. Rest of the proof is same. To show this, we use (\ref{mrinal}) in
place of (\ref{mandal}). By (\ref{mrinal}), it is enough to show that
$\ol \psi \ot R(X_m)$ can be lifted to a surjection from $\ol P\ot
R(X_m)\twoheadrightarrow \ol J\ot R(X_m)$. This is clearly true, since
$\ol J$ contains a monic polynomial in $X_m$ and $\ol P= \ol D\oplus P'$ by
(\ref{Wiemers}). \qed
\end{proof}
\medskip

The following lemma is very crucial to generalize above result.

\bl\label{semi} Let $R$ be a reduced ring and $D$ be a discrete Hodge
algebra over $R$. Let $L$ be a
rank one projective $D$-module. Then there exists a ring $S$ such that
\begin{enumerate}
 \item $R\hookrightarrow S\hookrightarrow Q(R)$,
\item $S$ is a finite $R$-module,
\item $R\hookrightarrow S$ is subintegral and
\item $L\ot_R S $ is extended from $S$.
\end{enumerate}
\el 

\proof Let $R\hookrightarrow B \hookrightarrow Q(R)$ be the
seminormalization of $R$. By Swan's result \cite[Theorem 6.1]{sw},
rank one projective modules over polynomial extensions of $B$ are
extended from $B$. Hence by (\ref{vorst}), rank one projective modules over
discrete Hodge algebras over $B$ are extended from $B$. In particular
$L\ot_R B$ is extended from $B$. By \cite[Theorem
  2.8]{sw}, $B$ is direct limit of $B_\lambda$, where
$R\hookrightarrow B_\lambda$ is finite and subintegral
extension. Since $L$ is finitely generated, we can find a subring
$S=B_\lambda$ for some $\lambda$ satisfying conditions $(1-4)$.
\qed
\medskip

We now prove the general case of (\ref{rem1}).

\bt\label{main0}
Let $R$ be a ring containing $\mathbb{Q}$ of dimension  $\geq 2$ and $D$ be a discrete Hodge algebra over $R$ of dimension $n> \text{dim}(R)$.
 Let $I$ be an ideal in $D$ of height $\geq 3$ and $P$ be any rank $n$ projective $D$-module.
 Suppose that there is a surjection $\alpha:P/IP\twoheadrightarrow I/I^2$. Then $\alpha$ can be lifted to a 
surjection $\beta: P \twoheadrightarrow I$.
\et
\proof

Without loss of generality, we may assume that $D$ is reduced. In particular, $R$ is reduced. 
Let $D=R[X_1,\cdots,X_m]/\mathcal{J}$, 
where $\mathcal{J}$ is an ideal of $R[X_1,\cdots,X_m]$ generated by square free monomials.
By (\ref{semi}), there exists an extension $R\hookrightarrow S$ such that 
\begin{enumerate}
 \item $R\hookrightarrow S\hookrightarrow Q(R)$,
\item $S$ is a finite $R$-module,
\item $R\hookrightarrow S$ is subintegral and
\item $\wedge^n(P)\ot_R S$ is extended from $S$.
\end{enumerate}
Let $E=S[X_1,\cdots,X_m]/\mathcal{J}$.  Since $\wedge^n(P)\ot_R S$ is
extended from $S$, by (\ref{rem1}), the induced surjection
$\alpha^*:P\ot E\twoheadrightarrow IE/I^2E$ can be lifted to a
surjection $\phi:P\ot E\twoheadrightarrow IE$.  By (\ref{Wiemers}),
$P=D\oplus Q$. In case $P=\wedge^n(P)\oplus D^{n-1}$, the rest of the proof
is given in \cite[Theorem 3.12]{dz1}. The proof of
\cite[Theorem 3.12]{dz1} works for $P=D\oplus Q$ also. Hence we are done.
\qed
\medskip

\section{Main Theorems: Codimension One Case:}
The aim of this section is to give an affirmative answer to Question
\ref{ques} mentioned in the introduction.  We start with the following
lemma which generalizes (\ref{dz2.13}).

\bl\label{lift} Let $R$ be a ring containing $\mathbb{Q}$ of dimension
$\geq 2$ and
$I$ be an ideal of $R[X,Y]$ of height $\geq 3$.
Let $P$ be a projective $R[X,Y]$-module of rank $\geq \dim (R)+1$ whose
determinant is extended from $R[X]$.  Suppose that there exists a
surjection $ \phi :P \twoheadrightarrow I/I^2$.  Then $\phi$ can be
lifted to a surjection $\psi : P \twoheadrightarrow I$.  \el

\proof If rank of $P$ is $>\dim (R)+1$, then we are done by (\ref{dz2.13}). 
So assume rank of $P$ $=\dim (R)+1$. Since $R$
contains $\mathbb{Q}$, using \cite[Lemma 3.3]{brs1} and replacing $Y$ by $Y-\lambda$ for some $\lambda\in \bq$, 
we can assume
that either $I(0)=R[X]$ or
$\text{ht}(I(0))=\text{ht}(I)$.
If $I(0)=R[X]$, then by (\ref{rem3.9}), we can lift $\phi$ to a surjection
$\alpha:P\twoheadrightarrow I/(I^2Y)$.

Now assume that $\text{ht}(I(0))=\text{ht}(I)\geq 3$. 
Let ``bar" denote the reduction modulo $Y$ and consider
 $\ol \phi:\ol P\twoheadrightarrow \ol I/{\ol I}^2$. 
By (\ref{dz2.13}), there exists a surjection
$\beta:\ol P\twoheadrightarrow \ol I$ which lifts
$\ol \phi$.
Therefore, again by (\ref{rem3.9}), we can lift $\phi$ to a 
surjection $\alpha:P\twoheadrightarrow I/(I^2Y)$. 
Therefore, in any case, we can lift  $\phi$ to a surjection $\alpha:P\twoheadrightarrow I/(I^2Y)$.

Consider the surjection $\alpha\ot R(Y): P\ot
R(Y)\twoheadrightarrow I\ot R(Y)/I^2\ot R(Y)$.
Since 
dim$(R(Y))=$ dim$(R)$,
by (\ref{dz2.13}), $\alpha\ot R(Y)$
can be lifted to a
surjection $\delta:P\ot R(Y)\twoheadrightarrow I\ot R(Y)$.  Using
(\ref{mrinal}), we get a surjection $\psi : P \twoheadrightarrow I$
which lifts $\alpha$ and hence lifts $\phi$.  \qed

\bp\label{main2} Let $R$ be a ring containing $\mathbb{Q}$ of
dimension $\geq 3$ and $D$ be a discrete Hodge algebra over $R$ of
dimension $d>\text{dim}(R)$. Let $I$ be an ideal in $D$ of height $\geq 4$
and $P$ be a projective $D$-module
of rank $n\geq max\{\text{dim}(R)+1, d-1\}$ whose determinant is extended from $R$. Suppose that
$\alpha:P\twoheadrightarrow I/I^2$ is a surjection.
Then there exists a surjection $\beta:P\twoheadrightarrow I$ that lifts $\alpha$.
\ep
\proof

 As in the proof of (\ref{main1}), we can assume that $R$ is reduced and
$D=R[X_1,\cdots,X_m]/\mathcal{I}$, where $\mathcal{I}$ is an ideal of
 $R[X_1,\cdots,X_m]$  generated by square free monomials. 
 where $X_{i_1}^{l_1}\cdots X_{i_k}^{l_k}\in \mathcal{I}$ and $l_i\geq
 1$ We prove the result using induction on $m$. If $m = 1$, then
 $D=R[X_1]$ and the result follows from (\ref{dz2.13}).

Let us assume that $m\geq 2$. If $D$ is a polynomial ring over $R$,
then we are done by (\ref{lift}).  Now suppose that $D$ is a
non-trivial discrete Hodge algebra. Then we can assume that $\mathcal
I=(\mathcal{K},X_m\mathcal{L})$, where $\mathcal{K}$ and
$\mathcal{L}$ are monomial ideals in $R[X_1,\ldots,X_{m-1}]$. 
Then $D=R[X_1,\cdots,X_m]/(\mathcal{K},X_m\mathcal{L})$.

Let $x_m$ and $L$ be the images of $X_m$ and $\mathcal{L}$ in $D$
respectively. We shall use ``tilde" when we move modulo $(x_m)$ and
``bar" when we move modulo $L$. We first go modulo $(x_m)$, i.e. to the
discrete Hodge algebra $\wi D=R[X_1,\cdots,X_{m-1}]/\mathcal{K}$ and
consider the surjection $\wi \alpha:\wi P\twoheadrightarrow
\iw/\iw^2$.  Note that $\widetilde I$ is an ideal of $\wi D$ of height
$\geq \dim (\wt D)-1$.  By induction hypothesis
on $m$, there exists a surjection $\phi:\wi P\twoheadrightarrow \iw$
which is a lift of $\wi\alpha$.  Therefore, using
(\ref{rem3.9}), we can lift $\alpha$ to a surjection
$\psi:P\twoheadrightarrow I/(I^2x_m)$.
 
We now move modulo $L$, i.e. $\ol D=\frac{R[X_1,\cdots,X_{m-1}]}{(\mathcal{K},\mathcal{L})}[X_m]:=D_0[X_m]$  and consider the surjection 
$$\ol \psi:\ol P\twoheadrightarrow \ol I/(\ol I^2X_m).$$

Observe that 
$\text{ht}(\ol I)\geq \text{dim}(R)\geq 3$.
If $\text{dim}(D_0)<n$, then by (\ref{dz2.13}), $\ol \psi$ can be lifted to a surjection $\theta:\ol P\twoheadrightarrow \io$. So assume
$\text{dim}(D_0)=n$. 
Since $\text{dim}(R(X_m))=\text{dim}(R)$ and $\ol D\ot R(X_m)=\frac{R(X_m)[X_1,\cdots, X_{m-1}]}{(\mathcal{K},\mathcal{L})}$,
by (\ref{rem1}), the surjection $\ol \psi \ot R(X_m) : \ol P\ot R(X_m) \twoheadrightarrow \ol I\ot R(X_m)/(\ol I^2 \ot R(X_m))$ can be lifted to
a surjection $\eta:\ol P\ot R(X_m)\twoheadrightarrow \ol I \ot R(X_m)$. 
By (\ref{mrinal}), there exists a surjection $\theta:\ol P\twoheadrightarrow \io$ which lifts $\ol\psi$.

Finally it follows from (\ref{rem3.9}) that there exists a surjection
$\beta:P\twoheadrightarrow I/(I^2x_mL)$ which
lifts $\psi$.
As $x_m L=0$ in $D$, we obtain a surjection 
$\beta:P\twoheadrightarrow I$ which lifts $\alpha$.
\qed
\medskip

Now we will answer Question \ref{ques0}. 

\bt\label{low} Let $R$ be a ring of dimension $\geq 3$ containing
$\mathbb{Q}$ and $D$ be a discrete Hodge algebra over $R$ of dimension
$d>\text{dim}(R)$. Let $I$ be an ideal in $D$ of height
$\geq 4$ and $P$ be a projective $D$-module of rank $n\geq
max\{\text{dim}(R)+1, d-1\}$. Suppose that $\alpha:P\twoheadrightarrow
I/I^2$ is a surjection. Then there exists a surjection
$\beta:P\twoheadrightarrow I$ which lifts $\alpha$.  
\et 

\proof
Without loss of
generality we may assume that $D$ is reduced. Using (\ref{move}), we
can lift $\alpha$ to a surjection $\alpha' :P\twoheadrightarrow I\cap
I_1$ such that $I+I_1=D$ and ht$(I_1)\geq n$. 

If ht$(I_1)> n$, then $I_1=D$ and hence $\alpha'$ is the required
surjective lift of $\alpha$. Assume ht$(I_1)= n$. The map $\alpha'$
induces a surjection $\alpha_1:P\twoheadrightarrow I_1/I_1^2$. If we can
show that $\alpha_1$ can be lifted to a surjection
$\Delta:P\twoheadrightarrow I_1$, then by subtraction principle
\cite[Proposition 3.2]{dk}, we can find a surjection
$\Delta_1:P\twoheadrightarrow I$ which lifts $\alpha$. Therefore it is
enough to show that $\alpha_1$ has a surjective lift $\Delta$. Now
replacing $I_1$ by $I$ and $\alpha_1$ by $\alpha$, we assume that
ht$(I)=n$.

 By (\ref{semi}), there exists an extension $R\hookrightarrow S$ such that
\begin{enumerate}
 \item $R\hookrightarrow S\hookrightarrow Q(R)$,
\item $S$ is a finite $R$-module,
\item $R\hookrightarrow S$ is subintegral and
\item $\wedge^n(P)\ot_R S$ is extended from $S$.
\end{enumerate}

Let $C$ be the conductor ideal of $R$ in $S$. Then $\text{ht}(C)\geq
1$.  Since $\hh(I)=n\geq max\{\text{dim}(R)+1, d-1\}$ and
$\hh(C)\geq1$, it follows that $\text{ht}(I^2\cap C)\geq 1$.
Therefore, we can choose an element $b\in I^2\cap C$ such that
$\text{ht}(b)=1$.  Let ``bar'' denote reduction modulo the ideal
$(b)$.  Consider the surjection $\overline\alpha : \ol P
\twoheadrightarrow \overline I/{\overline I}^2$ and note that
$\text{dim}(\ol R)< \text{dim}(R)$.

Now applying  (\ref{main0}), we can find a surjection 
$\gamma':\ol P\sur \ol{I}$ which lifts $\ol{\ga}$.
Choose a lift $\gamma:  P\lra I$ of $\gamma'$.
Since $b\in I^2$, $\gamma$ is a lift of $\alpha$ and hence $(\gamma(P),b)=I$.
Since hh$(I)=n$ and $b\in I^2$,
applying (\ref{EE}) and replacing $\gamma$ by $\gamma+b\delta$ for some
$\delta\in P^\ast$, we can assume that 
$\hh(\gamma(P))=n$.  

Applying (\ref{nak}), there exists an ideal $I^{\prime}$ of height
$\geq n$ such that $I'+ bD=D$ and $\gamma(P)=I\cap I'$. If
$\hh(I^{\prime})> n$, then $I'=D$ and hence $\gamma$ is the required
surjective lift of $\alpha$. Assume that ht$(I^{\prime})= n$ and
consider the surjection $\theta:P\twoheadrightarrow I'/{I'}^2$ induced
from $\gamma:P\twoheadrightarrow I\cap I'$.


Consider the surjection $\theta\ot_R S:P\ot S\twoheadrightarrow I'\ot
S/{I'^2\ot S}.$ Since $\wedge^n(P\ot_R S)$ is extended from $S$, by
(\ref{main2}), $\theta\ot S$ can be lifted to a surjection
$\Theta:P\ot S\twoheadrightarrow I'\ot S$. Now we need to show that we
get a surjection $\eta:P\twoheadrightarrow I'$ which lifts $\theta$.
In the case of $P=\wedge^n(P)\oplus D^{n-1}$, this is proved in
\cite[Lemma 5.1]{dz2}. Note that $P=D\oplus P'$, by
(\ref{Wiemers}). The proof of \cite[Lemma 5.1]{dz2} works in this case
also, so we do not repeat it here. Therefore we have a surjection
$\eta:P\twoheadrightarrow I'$ which lifts $\theta$.  Applying
subtraction principle \cite[Proposition 3.2]{dk}, we can find a
surjection $\beta:P\twoheadrightarrow I$ which lifts $\alpha$.  
\qed
\medskip

The following result is immediate from (\ref{low}).

\bc\label{cor}
Let $R$ be a ring of dimension $d\geq 3$ containing $\mathbb{Q}$ and $D=\frac{R[X_1,X_2,X_3]}{\mathcal{I}}$ be a discrete Hodge algebra over $R$. 
Let $I$ be an ideal in $D$ of height $\geq 4$
and $P$ be a projective $D$-module
 of rank $n\geq \text{dim}(R)+1$. Suppose that $\alpha:P\twoheadrightarrow I/I^2$ be a surjection. Then there exists a surjection $\beta:P\twoheadrightarrow I$ which 
 lifts $\alpha$.
\ec


\medskip

The following theorem is due to Ferrand and Szpiro \cite{sz}.

\bt\label{fs}
Let $R$ be a ring and $I\subset R$ be a locally complete intersection ideal of height $r\geq 2$ and $\text{dim}(R/I)\leq 1$.
Then there is a locally complete intersection ideal $J\subset R$ of height $r$ such that
\begin{enumerate}
 \item  $\sqrt I=\sqrt J$ and
\item $J/J^2$ is free $R/J$-module of rank $r$.
\end{enumerate}
\et

As an application of (\ref{low}), we improve a result of Mandal \cite[Corollary 2.2]{m2}, albeit with a stronger hypothesis on ideals.

\bt \label{set}
Let $R$ be a ring of dimension $\geq 3$ containing $\mathbb{Q}$ and $D$ be a discrete Hodge algebra over $R$ with 
$\text{dim}(D)=d>\text{dim}(R)$. Let $I$ be a locally complete intersection ideal in $D$ of height $n= max\{\text{dim}(R)+1, d-1\}$. 
Then there exist $f_1,\cdots,f_n\in I$ such that $\sqrt I=\sqrt {(f_1,\cdots,f_n)}$. In other words, $I$ is set 
theoretically generated by $n$ elements.
\et
\proof
By (\ref{fs}), there is a locally complete intersection ideal $J$ such that 
$\sqrt I=\sqrt J$ and
$J/J^2$ is a free $R/J$-module of rank $n$.
Applying (\ref{low}), we see that $J$ is generated by $n$ elements. Therefore, $I$ is set theoretically
generated by $n$ elements.
\qed
\medskip

\section{Some Auxiliary Results}

After answering Question \ref{ques} and Question \ref{ques0}, it is natural
to ask the following more general question. 

\quest
Let $R$ be a commutative Noetherian ring of dimension $\geq 1$ and 
$D$ be a discrete Hodge algebra over $R$ of dimension $n>\text{dim}(R)$.
 Let $I\subset D$ be an ideal of height $>\text{dim}(R)$. Suppose that 
$I=(f_1,\cdots,f_n)+I^2$, where $n\geq \text{dim}(D/I)+2$.
 Do there exist $g_1,\cdots,g_n\in I$ such that  $I=(g_1,\cdots,g_n)$ with $f_i-g_i\in I^2$? \smallskip

The above question has been answered affirmatively by Mandal when $D$ is a polynomial algebra over $R$ (\cite{m1}).
 Using \cite[Theorem 4.2]{drs} and following the proofs of (\ref{main1}) and (\ref{main2}), 
 we can obtain the following result which gives
 a partial answer to the above question.

\bt\label{main}
Let $R$ be a ring of dimension $d\geq 2$ containing $\mathbb{Q}$ and $D=\frac{R[X_1,\cdots,X_m]}{(J_1,X_mJ_2)}$,
  where $J_1,J_2$ are two ideals of $R[X_1,\cdots,X_{m-1}]$ generated by monomials. 
 Let $I$ be an ideal in $D$ of height $> d$. Suppose that $I=(f_1,\cdots,f_n)+I^2$ with $n\geq \text{dim}(D/I)+2$. Then $I=(g_1,\cdots,g_n)$ with $f_i-g_i\in I^2$ 
in each of the following cases:
\begin{enumerate}
 \item $n\geq max\{\text{dim}(D/J_1), \text{dim}(D/J_2)\}$ and $\text{ht}(\frac{I+J_2}{J_2})\geq 2$. 
\item $n= max\{\text{dim}(D/J_1)-1, \text{dim}(D/J_2)-1\}$ and $\text{ht}(\frac{I+J_2}{J_2})\geq 3$.  
\end{enumerate}
\et

As an application of (\ref{main}), we give some explicit examples. \smallskip

\example
Let $R$ be a ring of dimension $d\geq 4$ containing $\mathbb{Q}$ and $D=\frac{R[X_1,\cdots,X_4]}{(X_4J)}$ where $J=(X_1X_2, X_2X_3, X_1X_3)$. Let $I\subset D$ be an ideal of 
height $n\geq d+1$. Suppose that $I=(f_1,\cdots,f_n)+I^2$. Then there exist $g_1,\cdots g_n \in I$ such that $I=(g_1,\cdots g_n)$ with $f_i-g_i\in I^2$.
 In other words, the $n$-th Euler class group $E^n(D)$ is trivial. \smallskip

\proof
Using (\ref{main1}) and (\ref{main2}), we can assume that $n=d+1$. 
 We have $\text{dim}(D/J)=d+2$, i.e., $n=d+1=\text{dim}(D/J)-1$ and  $\text{ht}(\frac{I+J}{J})\geq3$.
Also note that $n=d+1\geq5 \geq \text{dim}(D/I)+2$. Now the result follows from (\ref{main}(2)).
\qed 
\smallskip

The following result follows from (\ref{main}). \smallskip

\example
Let $R$ be a ring of dimension $d\geq 3$ containing $\mathbb{Q}$ and $D=\frac{R[X_1,\cdots,X_m]}{(X_mJ)}$ where $J=(X_iX_j|1\leq i\neq j\leq m-1)$.
 Let $I\subset D$ be an ideal such that $\text{ht}(\frac{I+J}{J})\geq 3$. 
Suppose that $I=(f_1,\cdots,f_n)+I^2$ with $n\geq max\{d+1, \text{dim}(D/I)+2\}$. Then there exist $g_1,\cdots g_n \in I$ such that $I=(g_1,\cdots g_n)$ with $f_i-g_i\in I^2$.
\qed
\smallskip

Now using (\ref{cor}) and following the proof of (\ref{main}), we can
derive the following. \smallskip

\example
Let $R$ be a ring of dimension $d\geq 4$ containing $\mathbb{Q}$ and $D=\frac{R[X_1,\cdots,X_4]}{(J_1, X_4J_2)}$ 
where $J_1, J_2$ are two ideals in $R[X_1,X_2,X_3]$ generated by monomials and $\text{ht}(J_1+J_2)\geq 2$. Let $I\subset D$ be an ideal of 
height $n\geq d+1$. Suppose that $I=(f_1,\cdots,f_n)+I^2$. Then there exist $g_1,\cdots g_n \in I$ such that $I=(g_1,\cdots g_n)$ with $f_i-g_i\in I^2$.
 In other words, the $n$-th Euler class group $E^n(D)$ is trivial. \smallskip

\proof
Since $\text{dim}(D)\leq d+3$, the case $n\geq d+2$ is covered by (\ref{main1}) and (\ref{main2}). Let us assume that $n=d+1$. 
Then $n=d+1= \text{dim}(D/J_2)-1$ and $2n\geq \text{dim}(D)+2$. Now the result follows from (\ref{main}).
\qed

\smallskip
The following result follows from (\ref{main}). \smallskip

\example
Let $R$ be a ring of dimension $d\geq 4$ containing $\mathbb{Q}$ and $D=\frac{R[X_1,\cdots,X_5]}{(X_5J)}$ 
where $J=(X_1X_2X_3, X_1X_2X_4, X_2X_3X_4)$.
 Let $I\subset D$ be an ideal such that $\text{ht}(I)=n\geq d+1$. 
Suppose that $I=(f_1,\cdots,f_n)+I^2$ with $n\geq d+2$. Then there exist $g_1,\cdots g_n \in I$ such that $I=(g_1,\cdots g_n)$ with $f_i-g_i\in I^2$. In other words, the $n$-th Euler class group $E^n(D)$ is trivial.
\qed

\end{document}